\documentclass[11pt,a4paper]{article}
\usepackage{amsfonts,epsf,amsmath,amssymb,graphicx,tabularx,booktabs}
\usepackage[ruled]{algorithm2e}

\newtheorem{theorem}{\bf Theorem}[section]

\newtheorem{lemma}[theorem]{\bf Lemma}

\newcommand{\proof}{\noindent{\bf Proof.\ }}
\newcommand{\qed}{\hfill $\blacksquare$ \bigskip}

\begin{document}

\title{ \bf Generalizations of Wiener polarity index and terminal Wiener index}

\author{
Aleksandar Ili\' c \footnotemark[3] \\
Faculty of Sciences and Mathematics, Vi\v segradska 33, 18 000 Ni\v s \\
University of Ni\v s, Serbia \\
e-mail: \tt{aleksandari@gmail.com} \\
\and
Milovan Ili\' c \\
Faculty of Information Technology, Trg republike 3, 11 000 Beograd \\
University of Belgrade, Serbia \\
e-mail: \tt{ilic.milovan@gmail.com} \\
}

\date{\today}

\maketitle

\begin{abstract}
In theoretical chemistry, distance-based molecular structure descriptors are used for modeling
physical, pharmacologic, biological and other properties of chemical compounds. We introduce a
generalized Wiener polarity index $W_k (G)$ as the number of unordered pairs of vertices $\{u, v\}$
of $G$ such that the shortest distance $d (u, v)$ between $u$ and $v$ is $k$ (this is actually 
the $k$-th coefficient in the Wiener polynomial). For $k = 3$, we get
standard Wiener polarity index. Furthermore, we generalize the terminal Wiener index $TW_k (G)$ as
the sum of distances between all pairs of vertices of degree $k$. For $k = 1$, we get standard
terminal Wiener index. In this paper we describe a linear time algorithm for computing these
indices for trees and partial cubes, and characterize extremal trees maximizing the generalized
Wiener polarity index and generalized terminal Wiener index among all trees of given order $n$.
\end{abstract}

 \vskip 0.1cm

 \vskip 0.1cm

\section{Introduction}

Let $G = (V, E)$ be a connected simple graph with $n = |V|$ vertices and $m = |E|$ edges. For
vertices $u, v \in V$, the distance $d (u, v)$ is defined as the length of the shortest path
between $u$ and $v$ in $G$. The diameter $diam (G)$ is the greatest distance between two vertices
of $G$. Let $d_k (u)$ denotes the number of vertices on distance $k$ from the vertex $u$. Let $deg
(v)$ denotes the degree of the vertex $v$.

In theoretical chemistry molecular structure descriptors (also called topological indices) are used
for modeling physico-chemical, pharmacologic, toxicologic, biological and other properties of
chemical compounds \cite{GuPo88}. There exist several types of such indices, especially those based
on vertex and edge distances. Arguably the best known of these indices is the Wiener index $W$,
defined as the sum of distances between all pairs of vertices of the molecular
graph~\cite{DoEnGu01}
$$
W (G) = \sum_{u, v \in V (G)} d (u, v).
$$
Besides of use in chemistry, it was independently studied due to its relevance in social science,
architecture and graph theory. With considerable success in chemical graph theory, various
extensions and generalizations of the Wiener index are recently put forward
\cite{DaGuMuSw09,ToCo00}.

The Wiener polarity index of a graph $G$ is defined as the number of unordered pairs of vertices
$\{u, v\}$ of $G$ such that the shortest distance $d (u, v)$ between $u$ and $v$ is $3$,
$$
WP (G) = |\{ (u, v) \mid d (u, v) = 3, \ u, v \in V \}|.
$$
Hosoya \cite{Ho02} found a physico-chemical interpretation of $WP$. Du, Li and Shi \cite{DuLiSh09}
described a linear time algorithm for computing the Wiener polarity index of trees, and
characterized the trees maximizing the index among all trees of given order. Deng et
al.~\cite{De10,DeXi10,DeXiTa10} and Liu et al.~\cite{LiHoHu10} characterized extremal $n$-vertex
trees with given diameter, number of pendent vertices or maximum vertex degree.

For $k \geq 1$, we define the generalized Wiener polarity index as the number of unordered pairs of
vertices $\{u, v\}$ of $G$ such that the shortest distance $d (u, v)$ between $u$ and $v$ is $k$,
$$
W_k (G) = \frac{1}{2} \sum_{v \in V (G)} d_k (v) = |\{ (u, v) \mid d (u, v) = k, \ u, v \in V \}|.
$$

Notice that $W (G) = \sum_{k = 1}^{diam (G)} W_k (G)$. If $x$ is a parameter, then the Wiener polynomial of $G$ is defined as \cite{Ho88,SaYeZh96}
$$
W(G, x) = \sum_{u, v \in V (G)} x^{d (u, v)} = \sum_{k = 1}^{diam (G)} W_k (G) \cdot x^k.
$$

Therefore, the generalized Wiener polarity index is basically the $k$-th coefficient in the Wiener polynomial.

The terminal Wiener index of a graph $G$ is defined by Gutman, Furtula and Petrovi\' c in
\cite{GuFuPe09} as the sum of distances between all pairs of pendent vertices of $G$,
$$
TW (G) = \sum_{\substack{{u, v \in V (G)} \\ {\ deg (u) = deg (v) = 1}}} d (u, v).
$$
Furthermore, the authors described a simple method for computing TW of trees and characterized the
trees with minimum and maximum $TW$. Recently Deng and Zhang \cite{DeZh09} studied equiseparability on
terminal Wiener index. Independently, Sz\' ekely et al. \cite{SzWaWu10} introduced the same index
(the sum of distances between the leaves of a tree) and studied the correlation between various
distance-based topological indices.

For $k \geq 1$, we define the generalized terminal Wiener index as the sum of the distances between
all unordered pairs of vertices of degrees $k$,
$$
TW_k (G) = \sum_{\substack{{u, v \in V (G)} \\ {\ deg (u) = deg (v) = k}}} d (u, v).
$$

The paper is organized as follows. In Section 2 we introduce generalization of the Wiener polarity
index $W_k(G)$ and characterize the trees maximizing the generalized Wiener polarity index among all trees
of given order, while in Section 3 we designed linear algorithm for calculating this index. In
Section 4 we introduce generalization of the terminal Wiener index and characterize trees
maximizing the generalized terminal Wiener index among all trees of given order. In
Section 5 we present formula for calculation of $TW_k (G)$ for partial cubes and in particular
closed formula for $TW_3$ of coronene series $H_k$. We close the paper in Section 6 by proposing
new problems for research.

\section{Generalization of Wiener polarity index}

For $k = 1$, it can be easily seen that $W_1 (G) = m$, where $m$ is the number of edges. For $k =
2$, we have
$$
W_2 (G) = \sum_{v \in V} \binom{deg (v)}{2} = \frac{\sum_{v \in V} deg^2 (v)}{2} - m = \frac{M_1 (G)}{2} - m,
$$
where $M_1 (G)$ denotes the first Zagreb index of a graph \cite{NiKoMiTr03}.

For $k = 3$ we have the Wiener polarity index,
\begin{eqnarray*}
W_3 (T) &=& \sum_{uv \in E} (deg (v) - 1)(deg (u) - 1) = \sum_{uv \in E} deg (u) deg (v) - \sum_{v
\in V} deg^2 (v) + m \\
&=& M_2 (T) - M_1 (T) + m,
\end{eqnarray*}
where $M_2 (T)$ denotes the second Zagreb index of a graph \cite{IlSt10}.

In the following assume that $k \geq 3$. If the diameter of $T$ is less than $k$, then $W_k (T) =
0$. Therefore, the minimum value of $W_k (T)$ is zero, and it is achieved for all trees with $diam
(T) < k$ (for example the star $S_n$). On the other hand, we will prove that the maximum value of
$W_k (T)$ is achieved for a tree with diameter $k$ and with all pendent vertices on distance $k$.

The group of pendent vertices is defined as the set of all pendent vertices attached to the same
unique neighbor. Let $A_1$ and $A_2$ be two different groups of pendent vertices with the unique
neighbors $w_1$ and $w_2$, such that the distance between two arbitrary pendent vertices from these
groups is not equal to $k$. Let $p_1$ be the number of vertices on distance $k$ from an arbitrary pendent
vertex from $A_1$ and $p_2$ be the number of vertices on distance $k$ from an arbitrary pendent 
vertex from $A_2$. Without loss of generality assume that $p_1 \leq p_2$. 
If we remove all pendent vertices from $A_2$
and add them to the group $A_1$, we get a new tree $T'$ such that
$$
W_k (T') - W_k (T) = (|A_1|  p_1 +|A_2|  p_2) - (|A_1|  p_1 + |A_2| p_1) =|A_2| (p_2 - p_1) \geq 0.
$$
By repetitive application of this transformation, we will get a new tree with possibly increased
generalized Wiener polarity index. The diameter of $T'$ is not greater than the diameter of $T$
and each transformation introduces one new pendent vertex. By choosing two most distant
groups of pendent vertices, we will get the extremal tree with diameter equal to $k$. After that 
we can apply the transformation finitely many times, until all pendent vertices are on distance $k$ or $2$.

Assume that there are $p$ groups of pendent vertices with sizes $a_1, a_2, \ldots, a_p$ and $a_1 +
a_2 + \ldots + a_p = q$. Since $diam (T) = k$, we have $n - k + 1 \geq q \geq 2$. The distance
between any two pendent vertices not from the same group is equal to $k$, and therefore
$$
W_k (T) = \frac{1}{2} \sum_{i = 1}^p a_i (q - a_i) = \frac{1}{2} \left ( q^2 - \sum_{i = 1}^p
a_i^2 \right).
$$
The minimum value of $\sum_{i = 1}^p a_i^2$ under the condition $\sum_{i = 1}^p a_i = q$ is
achieved if and only if all numbers $a_i$ are as close as possible, i. e. $|a_i - a_j| \leq 1$ for
all $1 \leq i \leq j \leq p$. This can be easily proved by the transformation $(a_i, a_j) \mapsto
(a_i + 1, a_j - 1)$ with $a_j \geq a_i + 2$, since
$$
(a_i + 1)^2 + (a_j - 1)^2 - a_i^2 - a_j^2 = 2 (a_i - a_j) + 2 < 0.
$$

Notice that the tree is uniquely determined by the distances between pendent vertices \cite{Za65}.
A starlike tree is a tree with exactly one vertex of degree $\geq 3$. If $p = 2$, we have $W_k
(T) = a_1 a_2$ and $a_1 + a_2 = n - k + 1$ and finally
$$
W_k (T) = \left \lfloor \frac{n - k + 1}{2} \right \rfloor \cdot \left \lceil \frac{n - k + 1}{2}
\right \rceil
$$

Let $p > 2$. For $k$ odd, we can consider two groups of pendent vertices together with the unique
path connecting them. The third group of pendent vertices must be on equal distance from both
groups and that is impossible. Therefore, the extremal value for odd $k$ is achieved for $p = 2$.

For $k$ even, similarly it can be concluded that there is a unique tree with $p$ groups of pendent
vertices (starlike tree with $p$ paths with equal lengths together with groups of pendent vertices
attached at the end vertices of these $p$ paths). For $p > 2$, we have
$$
n = 1 + p \left (\frac{k}{2} - 1 \right) + \sum_{i = 1}^p a_i = 1 + p \left(\frac{k}{2} - 1 \right)
+ q,
$$
and since $p \leq q$, we have $p \leq 2 \cdot \frac{n - 1}{k}$. Therefore, using Cauchy--Schwartz
inequality it follows
\begin{eqnarray*}
W_k (T) &=& \frac{1}{2} \left ( q^2 - \sum_{i = 1}^p a_i^2 \right) \\
&\leq& \frac{1}{2} \left( q^2 - \frac{q^2}{p} \right ) \\
&=& \frac{1}{2} \left( n - 1 - \frac{pk}{2} + p \right)^2 \left (1 - \frac{1}{p} \right).
\end{eqnarray*}

Let
$$
f (p) = \frac{1}{2} \left( n - 1 - \frac{pk}{2} + p \right)^2 \left (1 - \frac{1}{p} \right),
$$
for $2 < p < 2 \cdot \frac{n - 1}{k} < 2 \cdot \frac{n - 1}{k - 2}$. The first derivative equals
$$
f' (p) = \frac{(2-2 n-2 p+k p) \left(2-2 n+2 p-k p-4 p^2+2 k p^2\right)}{8p^2}.
$$
Finally it holds that $f (p)$ is increasing for $2 < p < p^*$ and decreasing for $p^* \leq p <
2\cdot \frac{n -1}{k}$, where
$$
p^* = \frac{k - 2 + \sqrt{(k - 2)(16n + k - 18)}}{4(k - 2)} = \frac{1}{4} + \frac{1}{4}
\sqrt{\frac{16n + k - 18}{k - 2}}
$$
is the second largest root of $f' (p) = 0$. Therefore, the maximum of the Wiener polarity index for
$k$ even should be achieved for $p = 2$ or integers around $p^*$ (see Figure 1).

\begin{figure}[ht]
  \center
  \includegraphics [width = 3.5cm]{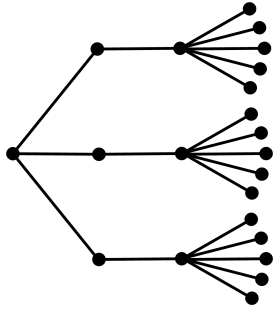}
  \caption { \textit{ The extremal tree with maximal generalized Wiener polarity index on $n = 22$ and $k = 6$. } }
\end{figure}

\section{Linear algorithm for $W_k$ of trees}

Let $T$ be an arbitrary tree rooted at the vertex $1$. We will process the vertices according to
the distance from the root vertex or in the recursive depth first search method \cite{CoLRS01}. For
each vertex $v$, keep the vector $a [v]$ of the length $k + 1$ that stores the number of
vertices in the subtree under $v$ on distances $0, 1, 2, \ldots, k$. It follows that $a [v][0] = 1$
and $a [v][1] = deg (v) - 1$ for all vertices different than the root. The matrix $a$ with
dimensions $n \times k$ is computed recursively in the first procedure.

In the second procedure we calculate the generalized Wiener polarity index. In DFS tree, for each
path of the length $k$ there is a unique vertex $v$ on the smallest distance from the root.
Therefore, we need to traverse all vertices $v$ and count the vertices that are on distance $k$ in
the subtree under $v$, such that the unique path connecting these vertices contains $v$. For the vertices
in the subtree under $v$, we just add $2 a [v][k]$. Otherwise, we need to consider all neighbors $u$ of
$v$ different than $parent [v]$ and for each $i = 0, 1, \ldots, k - 2$ count the number of vertex
pairs $(x, y)$ such that:
\begin{itemize}
\item $x$ is in the subtree under $u$ and $y$ is not;
\item $x$ is on the distance $i$ from $u$;
\item $y$ is on the distance $k - i - 1$ from $v$ and under $v$.
\end{itemize}
Finally we counted every vertex pair on distance $k$ twice, as showed in the second procedure. The time and memory
complexity is $O(nk)$. 

\begin{procedure}
\label{alg1} 
    \KwIn{The adjacency list of the tree $T$ with the root vertex $root$\,.}
    \KwOut{The array $parent$ and matrix $a$\,.}
    \BlankLine

    $a [v][0] = 1$\;
    \ForEach{neighbor $u$ of $v$}
    {
        \If{$(parent [u] = 0)$ and $(u \neq root)$}
        {
            $parent [u] = v$\;
            $DFS (u)$\;
            \For{$i = 0$ \KwTo $k - 1$}
            {
                $a [v][i + 1] = a [v][i + 1] +  a [u][i]$\;
            }
        }
    }
    \caption{DFS (vertex $v$)}
\end{procedure}

\begin{procedure}
\label{alg2} 
    \KwIn{The adjacency list of the tree $T$ with the root vertex $root$\, and the array $a$\,.}
    \KwOut{The generalized Wiener polarity index $Wk$\,.}
    \BlankLine

    $Wk = 0$\;
    \For{$v = 1$ \KwTo $n$}
    {
        $Wk = Wk + 2 \cdot a [v][k]$\;
        \ForEach{neighbor $u$ of $v$}
        {
            \If{$parent [v] \neq u$}
            {
                \For{$i = 0$ \KwTo $k - 2$}
                {
                    $Wk = Wk + a [u][i] \cdot (a [v][k - 1 - i] - a [u][k - 2 - i])$\;
                }
            }
        }
    }
    \Return{$Wk / 2$}\;
    \caption{GWP (G)}
\end{procedure}

\section{Generalization of terminal Wiener index}

For $k$ regular graphs, $TW_k (G) = W (G)$ and $TW_i (G) = 0$ for $i \neq k$. For $k = 1$, we have
terminal Wiener index.

\begin{theorem}[\cite{DoEnGu01}]
\label{th:wiener} Let $T$ be a tree on $n$ vertices. Then,
$$
(n - 1)^2 = W (S_n) \leq W (T) \leq W (P_n) = \binom{n+1}{3},
$$
with equality if and only if $T \cong S_n$ or $T \cong P_n$.
\end{theorem}

\begin{theorem}[\cite{IlIlSt10}]
\label{thm-pi} Let $w$ be a vertex of a nontrivial connected graph $G$. For nonnegative integers
$p$ and $q$, let $G (p, q)$ denote the graph obtained from $G$ by attaching to vertex $w$ pendent
paths $P = w v_1 v_2 \ldots v_p$ and $Q = w u_1 u_2 \dots u_q$ of lengths $p$ and $q$,
respectively. If $p \geq q \geq 1$, then
$$
W (G (p, q)) < W (G (p + 1, q - 1)).
$$
\end{theorem}

For $k = 1$, Gutman, Furtula and Petrovi\' c in \cite{GuFuPe09} characterized the extremal trees
within the class of all trees with $n$ vertices that maximize and minimize terminal Wiener index.
The path $P_n$ is the unique tree that minimizes $TW_1$. 

For $k = 2$, the unique tree that maximizes $TW_2$ is the path $P_n$. 
This follows from Theorems \ref{th:wiener} and the simple fact that in every tree there are at least 2 pendent vertices
(there are at most $n - 2$ vertices of degree two, and the maximum sum between all such pairs is achieved for path).

Let $T^*$ be the extremal tree that maximizes the generalized terminal Wiener index for $k \geq 3$.
It can be easily proved that there are no pendent paths $P = v_1 v_2 \ldots v_p$ of length greater
than 2or equal to attached at some vertex of $T^*$. Otherwise, remove pendent edges one by one and subdivide
any edge $e$ such that both components of $T^* - e$ contain vertices of degree $k$ -- this way we
increase $TW_k$. If such edge does not exist, remove the edge $v_{p-1}v_p$ and add new edge
$v_{p-2} v_p$ -- this way the generalized terminal Wiener index remains the same, or increases if
and only if $k = 3$.

\begin{lemma}
Among trees on $n$ vertices the maximal possible number of vertices of degree $k \geq 2$ is
$$m (n, k) = \left \lfloor \frac{n - 2}{k - 1} \right \rfloor.$$
\end{lemma}

\proof Consider the induced graph $H$ composed of the vertices of degree $k$. Let $h$ be the number
of vertices in $H$, and let $f$ be the number of edges in $H$. Since $H$ is acyclic and possible
disconnected, we have $0 \leq f \leq h - 1$. Furthermore, for the number of edges in tree $T$ holds
$$
h \cdot k - (h - 1) \leq h \cdot k - f \leq n - 1.
$$
It follows that $h (k - 1) \leq n - 2$, which completes the proof. \qed

\begin{figure}[ht]
  \center
  \includegraphics [width = 10cm]{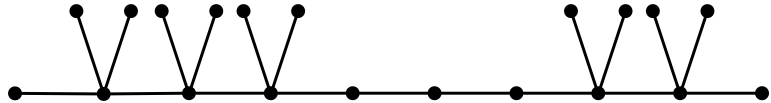}
  \caption { \textit{ The extremal tree with maximal generalized terminal Wiener index with $n = 20$ and $k = 4$ ($s = 8$ and $p = 5$) } }
\end{figure}

Let $C_{n, k, p}$ be the caterpillar obtained from a path of length $s + 2 = n - p \cdot (k - 2)$,
by attaching $k - 2$ pendent vertices to exactly $p$ vertices of a path $P_{s + 2} = v_0 v_1 \ldots
v_{s} v_{s + 1}$, starting from the vertices $v_1$ and $v_{s}$ on the both ends towards center
symmetrically (see Figure 2). The generalized terminal Wiener index can be easily calculated by 
considering the pairs of vertices $(v_1, v_s), (v_2, v_{s-1}), ...$ and summing the distance between all intermediate vertices:

\begin{eqnarray}
TW_k (C_{n, k, p}) &=& \sum_{1 \leq i, j, \leq s} d (v_i, v_j) \nonumber \\
&=& (p - 1)(s - 1) + (p - 3)(s - 3) + (p - 5)(s - 5) + \ldots \nonumber \\
&=& \left\{
\begin{array}{l l}
  \frac{1}{12} p (3ps - p^2 - 2), & \qquad \mbox{if $p$ is even} \\
  \frac{1}{12} (p + 1)(p - 1)(3s - p), & \qquad \mbox{if $p$ is odd.}
\end{array} \right. \nonumber \\
&=& \left\{
\begin{array}{l l}
  \frac{1}{12} p \left(3 n p + 5 p^2 - 3 k p^2 - 2 - 6 p \right), & \qquad \mbox{if $p$ is even} \\
  \frac{1}{12} (p + 1)(p - 1)(3n + 5p - 3kp - 6), & \qquad \mbox{if $p$ is odd.}
\end{array} \right. \label{eq:c-formula}
\end{eqnarray}

\begin{figure}[ht]
  \center
  \includegraphics [width = 9cm]{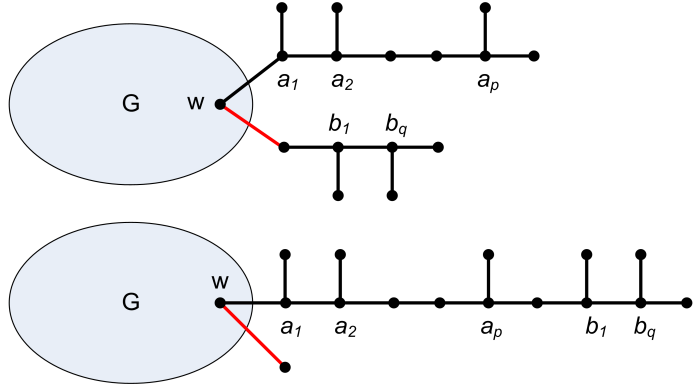}
  \caption { \textit{ Transformation that increases generalized terminal Wiener index $TW_3$. } }
\end{figure}

We call a caterpillar $C$ 3-bounded if all vertices of $C$ have degree less than or equal to $3$.

\begin{theorem}
Let $T$ be a tree on $n > 4$ vertices. Then
$$
TW_3 (T) \leq TW_3 (C_{n, 3, \lfloor n / 2 \rfloor - 1}),
$$
with equality if and only if $T \cong C_{n, 3, \lfloor n / 2 \rfloor - 1}$.
\end{theorem}

\proof Let $T^*$ be a rooted tree with maximal value of generalized terminal Wiener index for $k =
3$. We can also assume that the number of vertices of degree 3 is greater than two.

If $T^*$ is not 3-bounded caterpillar, consider a branching vertex $w$ such that in the subtree under
$w$ there are only 3-bounded caterpillars attached at $w$ 
(it can happen that there are only pendent vertices attached to $w$). 
Let $C_1$ and $C_2$ be two caterpillars
attached at $w$, such that $C_1$ has $p$ vertices $v_1, v_2, \ldots, v_p$ of degree 3 and $C_2$ has
$q$ vertices $u_1, u_2, \ldots, u_q$ of degree $3$ (see Figure 3). Without loss of generality, we
can assume that the number of vertices of degree $3$ in $C_1$ is greater than or equal to the number 
of vertices of degree $3$, namely $p \geq q$. 

Let $a_1, a_2, \ldots, a_p$ be the distances from the vertex $w$ to the vertices $v_1, v_2, \ldots,
v_p$ and $b_1, b_2, \ldots, b_q$ be the distance from the vertex $w$ to the vertices $u_1, u_2,
\ldots, u_q$. Let $D (w)$ be the sum of distances from the vertex $w$ to the vertices of degree $3$
in the subgraph $G$ (see Figure~3).

The generalized terminal Wiener index of the tree $T^*$ equals
\begin{eqnarray*}
TW_3 (T^*) &=& TW_3 (G) + (p + q) D (w) + r \left ( \sum_{i = 1}^p a_i + \sum_{j = 1}^q b_j \right)
\\ && + \sum_{i < j} (a_j - a_i) + \sum_{i < j} (b_j - b_i) + q \cdot \sum_{i = 1}^p a_i + p \cdot \sum_{j = 1}^q b_j,
\end{eqnarray*}
where $r \geq 1$ is the number of vertices of degree 3 in $G$.

After reattaching the caterpillar $C_2$ to the end of caterpillar $C_1$, the degree of vertex $w$
remains the same. The generalized terminal Wiener
index of transformed tree $T'$ equals
\begin{eqnarray*}
TW_3 (T') &=& TW_3 (G) + (p + q) D (w) + r \cdot \sum_{i = 1}^p a_i + r \left ( q \cdot a_p +
\sum_{j = 1}^q b_j  \right)
\\ && + \sum_{i < j} (a_j - a_i) + \sum_{i < j} (b_j - b_i) + q \cdot \left ( p \cdot a_p - \sum_{i = 1}^p a_i \right) + p \cdot \sum_{j = 1}^q b_j
\end{eqnarray*}

By subtraction, we get
$$
TW_3 (T') - TW_3 (T^*) = r q \cdot a_p + p q \cdot a_p - 2 q \cdot \sum_{i = 1}^p a_i,
$$
Since $a_p \geq a_i$ for $i = 1, 2, \ldots, p$, we have
\begin{eqnarray*}
TW_3 (T') - TW_3 (T^*) &=& q \left( a_p (r + p) - 2 \sum_{i = 1}^p a_i \right) \\
&\geq& 2q \cdot \sum_{i = 1}^p (a_p - a_i) \geq 0,
\end{eqnarray*}
with equality if and only if $p + r = 2$. Since $p \geq q \geq 1$ and $r \geq 1$, the equality holds iff $p = q = r = 1$.

It follows that using this transformation we can not decrease the generalized terminal Wiener index $TW_3$.
If $deg (w) > 3$, we can remove the pendent edge and subdivide some edge $e$ such that both
components of $T^* - e$ contain vertices of degree $3$, and further increase the generalized Wiener
index. Therefore after these transformations, there is only one 3-bounded caterpillar under $w$.

Finally, we conclude that the extremal tree is a 3-bounded caterpillar, and it can be easily seen that the maximal value of
$TW_3$ is achieved for caterpillar of the form $C_{n, k, p}$. 

For $k = 3$, from \eqref{eq:c-formula} we have
$$
f (p) = \left\{
\begin{array}{l l}
  \frac{1}{12} p \left(3 n p - 4 p^2 - 2 - 6 p \right), & \qquad \mbox{if $p$ is even} \\
  \frac{1}{12} (p + 1)(p - 1)(3n - 4p - 6), & \qquad \mbox{if $p$ is odd.}
\end{array} \right.
$$
and
$$
f (p) - f (p - 2) = \left\{
\begin{array}{l l}
  (p - 1)(n - 2p) - 1, & \qquad \mbox{if $p$ is even} \\
  (p - 1)(n - 2p), & \qquad \mbox{if $p$ is odd,}
\end{array} \right.
$$
which is greater than zero since $n \geq 2p + 2$. By direct verification, we get $f (\lfloor
\frac{n}{2} \rfloor - 1) > f (\lfloor \frac{n}{2} \rfloor - 2)$ for $n > 4$.

Therefore, it holds
$$
TW_3 (T) \leq TW_3 (C_{n, 3, \lfloor n / 2 \rfloor - 1}),
$$
with equality if and only if $T \cong C_{n, 3, \lfloor n / 2 \rfloor - 1}$. \qed

This can be further generalized -- for all $k > 3$ the caterpillar that maximizes the function $f (p) = TW_k (C_{n,
k, p})$ has maximal generalized terminal Wiener index among trees on $n$ vertices.

\section{Calculating $TW_k$ of partial cubes}

The $n$-cube $Q_n$ is the graph whose vertex set consists of all
binary $n$-tuples (hence $|V(Q_n)|=2^n$), two vertices are
adjacent if the corresponding tuples differ in precisely
one position. The central metric feature of the $n$ cube is the
fact that the distance between two vertices is equal to the
number of positions in which they differ. 
A subgraph $H$ of a graph $G$ is called {\em isometric}
if for any vertices $u$ and $v$ of $H$, $d_H(u,v) = d_G(u,v)$.
{\em Partial cubes} are isometric subgraphs of hypercubes.
Important examples of partial cubes are hypercubes, even cycles,
(chemical) trees, median graphs, benzenoid systems, phenylenes.
The Cartesian product of partial cubes is again a partial cubes.

Let $G$ be a connected graph. Then $e = xy$ and $f = uv$ are in the Djokovi\' c--Winkler $\Theta$
relation \cite{Wi84} if
$$
d (x, u) + d (y, v) \neq d (x, v) + d (y, u).
$$

The relation $\Theta$ is always reflexive and symmetric, and is transitive on partial cubes.
Therefore, $\Theta$ partitions the edge set of a partial cube $G$ into equivalence classes $F_1,
F_2, \ldots, F_s$, called $\Theta$-classes (or cuts). For any $1 \leq i \leq s$, the graph $G-F_i$ consists
of two connected components. The vertex sets of these components will be denoted with $W_{(i,0)}$ and
$W_{(i,1)}$, because they can be described as the vertices whose $i$-th coordinate is 0 and 1,
respectively. The sets $W_{(i,\chi)}$, $1\leq i\leq s$, $\chi\in \{0,1\}$ are called halfspaces of
$G$, while $W_{(i,0)}$ and $W_{(i,1)}$ are complementary halfspaces \cite{ImKl00} and it holds
$|W_{(i,0)}| + |W_{(i,1)}| = n$.

Let $G$ be a partial cube with halfspaces $W_{(i,\chi)}$, $1\leq i\leq s$, $\chi\in \{0,1\}$. For any
$1\leq i\leq s$ and any $\chi\in \{0,1\}$, $|W_{(i,\chi)}|^{(k)}$ denotes the number of vertices
with degree $k$ in $W_{(i, \chi)}$. 

\begin{theorem}
\label{thm:twk} Let $G$ be a partial cube with halfspaces $W_{(i,\chi)}$, $1\le i\le s$, $\chi\in
\{0,1\}$. Then
$$
TW_k (G) = \sum_{i=1}^s |W_{(i,0)}|^{(k)} \cdot |W_{(i,1)}|^{(k)}.$$
\end{theorem}

\proof $G$ is a partial cube, hence vertices of $G$ can be considered as a binary $s$-tuple
$u=u_1u_2\ldots u_s$. Moreover, since $G$ is isometric in $Q_s$, the distance between two vertices
is the number of positions in which they differ. Set $\delta(x,y) = 0$ if $x=y$, and $\delta(x,y) =
1$ for $x\not=y$. Then
\begin{eqnarray*}
TW_k (G) & = & \sum_{\substack{{u, v \in V (G)} \\ {\ deg (u) = deg (v) = k}}} d (u, v) \\
      & = &\sum_{\substack{{u, v \in V (G)} \\ {\ deg (u) = deg (v) = k}}} \sum_{i=1}^s \delta (u_i, v_i)\\
      & = & \sum_{i=1}^s \left(  \sum_{\substack{{u, v \in V (G)} \\ {\ deg (u) = deg (v) = k}}} \delta (u_i, v_i) \right) \\
      & = & \sum_{i=1}^s |W_{(i,0)}|^{(k)} \cdot |W_{(i,1)}|^{(k)}\, .
\end{eqnarray*}

This completes the proof. \qed

An advantage of Theorem~\ref{thm:twk} comparing to computing $TW_k(G)$ by the definition is that we
do not need to compute distances, but only to count vertices in the classes. This theorem can be
considered as another instance of Klav\v zar "cut method". For its general description and an
overview of its applications in chemical graph theory see survey \cite{Kl08}.

\begin{figure}[ht]
  \center
  \includegraphics [width = 8.5cm]{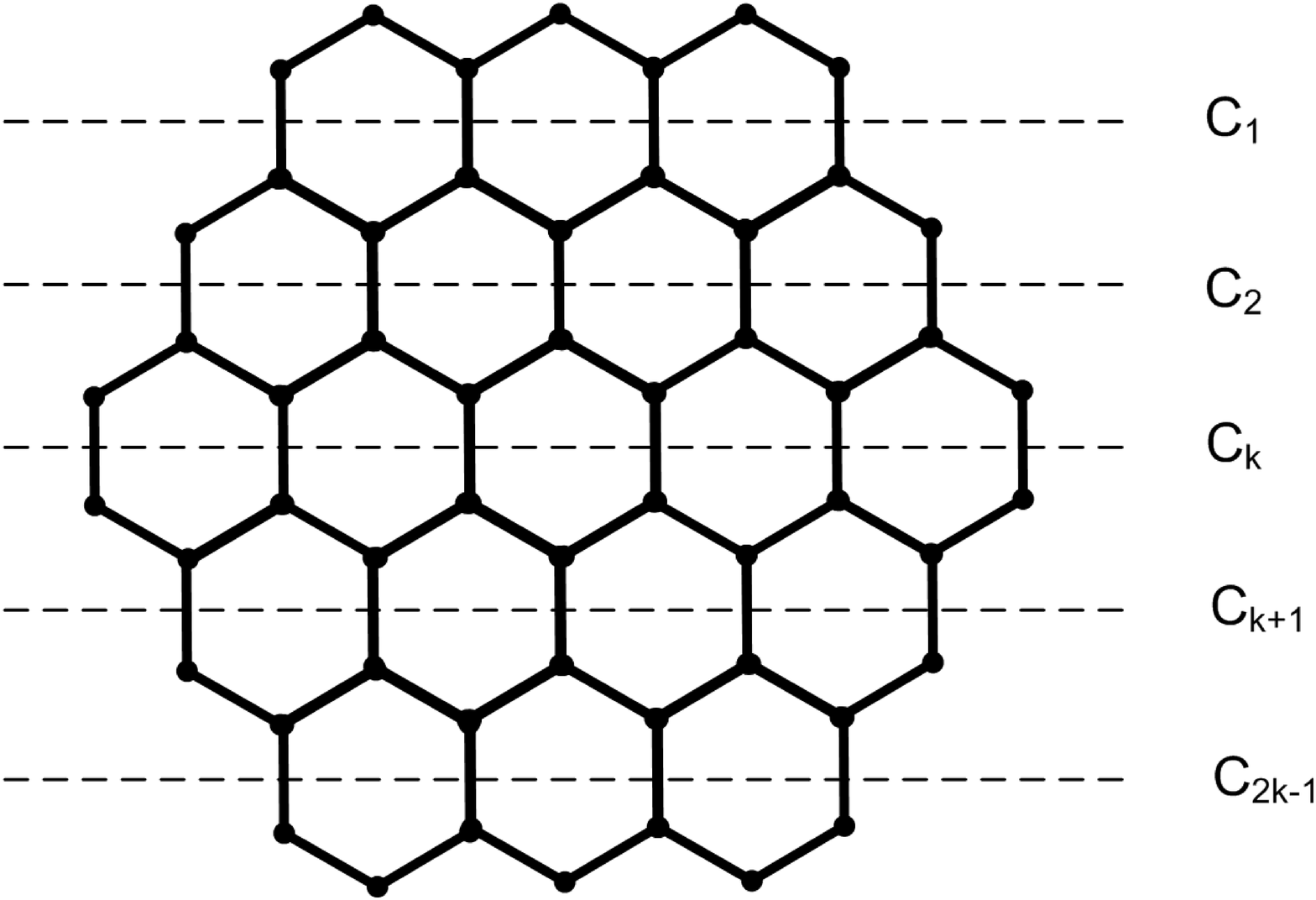}
  \caption{The Coronene / Circumcoronene $H_3$.}
\label{fig:coronene}
\end{figure}

As an example, we obtain a closed expression for $TW_3$ of the coronene/circumcoronene homologous
series $H_k$. In Figure \ref{fig:coronene}, $2k - 1$ horizontal elementary cuts of $H_k$ are
presented. There exist two additional groups of $2k - 1$ equivalent cuts, obtained by rotating the
former group by $\frac{\pi}{3}$ and $-\frac{\pi}{3}$. The number of vertices of $H_k$ equals $n_k =
6k^2$, while there are exactly $6k$ vertices of degree two.

Using symmetry, the contribution of the elementary cut $C_i$ is equal to the contribution of $C_{2k
- i}$, $i = 1, 2, \ldots, k - 1$. By induction it follows that for $i = 1, 2, \ldots, k$ the number
of vertices above cut $C_i$ equals $i (2k + i)$, while the number of vertices of degree 2 equals $k
+ 2i$. Therefore by using Theorem \ref{thm:twk} and described cuts, we have
\begin{eqnarray*}
\frac{1}{3} TW_3 (H_k) &=& (3k^2 - 3k)^2 + 2 \sum_{i = 1}^{k - 1} (2ki + i^2 - k - 2i)(6k^2 - 6k - 2ki - i^2 + k + 2i)\\
                       &=& \frac{164k^5}{15} - \frac{82 k^4}{3} +\frac{58 k^3}{3} -\frac{5 k^2}{3} -\frac{19k}{15}.
\end{eqnarray*}
Finally, we derive the fifth-order polynomial formula for the generalized terminal Wiener index of
$H_k$
$$
TW_3 (H_k) = \frac{1}{5} (k - 1) k (2k - 1)( 82k^2 - 82k - 19).
$$

Using similar methods as in \cite{ChKl97,IlKlSt10,KlGu97} one can obtain the closed formulas for
other chemical graphs (trees, benzenoid chains, phenylenes, ...) and design a linear algorithm for
$TW_k$ of benzenoid systems.

\section{Concluding remarks}

The Wiener polarity index and the terminal Wiener index are very new molecular-structure
descriptors and only a limited number of mathematical and chemical properties were established so
far. In this paper we generalized these indices and open new perspectives for the future research.

Another generalization of these indices may be the following
$$
W_k^* (G) = |\{ (u, v) \mid d (u, v) \leq k, \ u, v \in V \}| = W_1 (G) + W_2 (G) + \ldots + W_k
(G)
$$
and
$$
TW_k^* (G) = \sum_{\substack{{u, v \in V (G)} \\ {\ deg (u) \leq k, \ deg (v) \leq k}}} d (u, v)
\geq TW_1 (G) + TW_2 (G) + \ldots + TW_k (G).
$$

It would be nice to study mathematical and algorithmic properties of these indices and report their
chemical relevance. These indices are obtained from the famous Wiener index, which has many
applications in chemistry, graph theory and computer science.

\bigskip {\bf Acknowledgement. } This work was supported by Research Grants 174010 and
174033 of Serbian Ministry of Science. The authors are grateful to Nikola Milosavljevi\' c, Emeric Deutsch and Ivan Gutman for 
several useful suggestions while preparing the article.

\end{document}